\documentclass[11pt]{amsart}
\usepackage{verbatim}
\usepackage{amssymb}

\newtheorem{thm}{Theorem}[section]
\newtheorem{prop}[thm]{Proposition}
\newtheorem{lem}[thm]{Lemma}

\theoremstyle{definition}

\theoremstyle{remark}
\begin{document}
\title[Locally nilpotent varieties of groups]{Certain locally nilpotent varieties of groups}
\author{ Alireza Abdollahi}
  \address{Department of Mathematics\\ University of Isfahan;Isfahan 81746-73441;
  Iran; and
   Institute for Studies in Theoretical Physics and Mathematics, Tehran, Iran.}
  \email{a.abdollahi@sci.ui.ac.ir}
  \thanks{This research was in part supported by a grant from IPM}
\maketitle
 \begin{abstract}
 Let $c\geq 0$,  $d\geq 2$ be integers and $\mathcal{N}_c^{(d)}$ be the variety of groups in which every
$d$-generator subgroup is nilpotent of class at most $c$. N.D. Gupta posed this question that for what values
of $c$ and $d$ it is true that $\mathcal{N}_c^{(d)}$ is locally nilpotent?
We prove that if $c\leq 2^d+2^{d-1}-3$ then
the variety $\mathcal{N}_c^{(d)}$ is locally nilpotent and  we reduce  the question of Gupta about  the
periodic groups in $\mathcal{N}_c^{(d)}$ to the   prime power finite exponent groups in this   variety.
\end{abstract}
\maketitle
\theoremstyle{plain}
\newtheorem*{Prov}{Provision}
\newtheorem {mythm}{Theorem}[section]
\newtheorem {mylem}{Lemma}
\section{Introduction and results}
Let $c\geq 0$,  $d\geq 2$ be integers and  $\mathcal{N}_c$ be the variety of nilpotent groups of
class at most
$c$. We denote by $\mathcal{N}_c^{(d)}$ the variety of groups in which every $d$-generator subgroup
is in $\mathcal{N}_c$. In \cite{Gupta}, Gupta posed the following question:\\
 For what values
of $c$ and $d$ it is true that $\mathcal{N}_c^{(d)}$ is locally nilpotent?\\
Then he proved that for $c\leq (d^2+2d-3)/4$, the variety $\mathcal{N}_c^{(d)}$ is locally nilpotent.
In  \cite{Endimioni},  Endimioni improved the latter result where he proved that for $c\leq 2^d-2$,
the variety $\mathcal{N}_c^{(d)}$ is locally nilpotent.
Here we improve the number $2^d-2$ to $2^d+2^{d-1}-3$. In fact we prove:
\begin{thm} \label{thm0} ~~\\
(1) \; For $c\leq 2^d+2^{d-1}-3$, the variety $\mathcal{N}_c^{(d)}$ is locally nilpotent.\\
(2) \; For $c\leq 2^d+2^{d-1}+2^{d-2}-3$, every $p$-group in the  variety $\mathcal{N}_c^{(d)}$
is locally nilpotent, where $p\in\{2,3,5\}$.
\end{thm}
Note that  the variety $\mathcal{N}_c^{(2)}$ is contained in the variety of  $c$-Engel groups and it is yet
unknown whether every $c$-Engel group is locally nilpotent, even,  so far there is no published example
of a  non-locally nilpotent group in the variety $\mathcal{N}_c^{(2)}$.
 In the last section of this paper, by considering the
problem of locally nilpotency of the variety $\mathcal{N}_c^{(2)}$, we study  periodic groups
in this variety. Note that since every two generator subgroup of a group in $\mathcal{N}_c^{(2)}$ is
nilpotent, every periodic group in $\mathcal{N}_c^{(2)}$ is a direct product of $p$-groups ($p$  prime).
 We reduce the question of Gupta for  periodic
groups  in  $\mathcal{N}_c^{(d)}$ to the locally nilpotency of $p$-groups of finite exponent in
 this variety where the exponent depends only on the  numbers $p$  and $c$. In fact we prove that
\begin{thm}\label{thm1}
Let $p$ be a prime, $c>1$ an integer and $r=r(c,p)$ be the integer such that $p^{r-1}<c-1\leq p^r$.
Then the following are equivalent:\\
(1)\; every $p$-group in  $\mathcal{N}_c^{(2)}$ is locally nilpotent.\\
(2)\; if $p$ is odd,  every $p$-group of exponent dividing $p^r$ in $\mathcal{N}_c^{(2)}$ is
locally nilpotent, and if  $p=2$, every
$2$-group  of exponent dividing $2^{r+1}$ in  $\mathcal{N}_c^{(2)}$  is  locally nilpotent.
\end{thm}
\section{Groups in the variety $\mathcal{N}_c^{(d)}$}
Let $F_{\infty}$ be the free group of infinite countable rank on the set  $\{x_1,x_2,\dots\}$, we define
inductively  the following words  in $F_{\infty}$:\\
$$W_1=W_1(x_1,x_2)=[x_1,x_2,x_1,x_2],$$
$$W_n=W_n(x_1,x_2,\dots,x_{n+1})=[W_{n-1},x_{n+1},W_{n-1},x_{n+1}] \;\;\;\; n>1;$$
$$V_1=V_1(x_1,x_2,x_3)=[[x_2,x_1,x_1,x_1,x_1],x_3,[x_2,x_1,x_1,x_1,x_1],x_3],$$
$$V_n=V_n(x_1,x_2,x_3,\dots,x_{n+2})=[V_{n-1},x_{n+2},V_{n-1},x_{n+2}] \;\;\;\; n>1.$$
For a group $G$ and a subgroup  $H$ of $G$, we denote by $HP(G)$ the Hirsch-Plotkin radical of $G$ and
$H^G$ the normal closure of $H$ in $G$.
We use the following result due to Heineken (see Lemma 8 of \cite{Heineken} and see Lemma 2 of
\cite{newell} for the left-normed version).
\begin{lem} \label{thmhein}
Let $G$ be a group and $g$ an element in $G$ such that $[g,x,g,x]=1$ for all $x\in G$.  Then the normal
closure of $\left<g\right>$ in $\left<g\right>^G$ is abelian. In particular, $g\in HP(G)$.
\end{lem}
\begin{lem}\label{lem00}
Let $G$ be a group satisfying the law $W_n=1$ for some integer $n\geq 1$. Then, $G$ has a normal series
$1=G_n\lhd G_{n-1}\lhd\cdots\lhd G_1=G$ in which each factor $G_i/G_{i+1}$ is locally nilpotent
($i=1,2,\dots,n-1$).
\end{lem}
\begin{proof}
We argue by induction on $n$. If $n=1$, then Lemma \ref{thmhein}
yields that $x_1\in HP(G)$ for all $x_1\in G$ and so $G$ is
locally nilpotent. Now suppose that the lemma is true for $n$ and
 $G$ satisfies the law $W_{n+1}=1$.  By Lemma \ref{thmhein}, we have
 $W_n(x_1,\dots,x_{n+1})\in HP(G)$ for all $x_1,\dots,x_{n+1}\in G$
and so $G/HP(G)$ satisfies the law $W_n=1$. Thus by induction hypothesis $G/HP(G)$ has a normal series
 of length $n$  with locally nilpotent factors, it completes the proof.
\end{proof}
\begin{lem} \label{lem000}
Let $G$ be a $p$-group satisfying the law $V_n=1$ for some integer $n\geq 1$ where $p\in\{2,3,5\}$. Thus
$G$ has a normal series $1=G_{n+1}\lhd G_n\lhd\cdots\lhd G_1=G$ in which each factor $G_i/G_{i+1}$ is
 locally nilpotent ($i=1,2,\dots,n$).
\end{lem}
\begin{proof}
We argue by induction on $n$. If $n=1$, then by Lemma
\ref{thmhein} \\ $[x_2,x_1,x_1,x_1,x_1] \in HP(G)$ for all
$x_1,x_2\in G$. Thus $G/HP(G)$ is a $4$-Engel group and since
every $4$-Engel $p$-group is locally nilpotent where
$p\in\{2,3,5\}$ (see Traustason \cite{Tr} and Vaughan-Lee
\cite{VauLee}), so $G/HP(G)$ is locally nilpotent. Now suppose
that $G$ satisfies the law $V_{n+1}=1$. By Lemma \ref{thmhein},
$V_n(x_1,\dots,x_{n+2})\in HP(G)$ for all $x_1,\dots,x_{n+2}\in
G$,  so $G/HP(G)$ satisfies the law $V_n=1$. Thus by induction
hypothesis, $G/HP(G)$ has a normal series of length $n+1$ with
locally nilpotent factors, it completes the proof.
\end{proof}
We use  in the sequel  the following special case  of this well-known  fact  due
 to Plotkin \cite{Plot}
that every Engel radical group is locally nilpotent. (see also Lemma 2.2 of \cite{Endimioni})
\begin{lem} \label{remark}
Let $H$ be a normal subgroup of an Engel group $G$. If $H$ and $G/H$ are locally nilpotent,
 then $G$ is locally nilpotent.
\end{lem}

{\bf Proof of Theorem \ref{thm0}.}
 One can see that $W_n$ and $V_n$ are, respectively,  in  the  $(2^{n+1}+2^n-2)$th term  and
 $(2^{n+2}+2^{n+1}+2^n-2)$th term of the lower central series of $F_{\infty}$,  for
all integers $n\geq 1$.\\ (1) \; every group $G$ in the variety
$\mathcal{N}_c^{(d)}$ satisfies the law $W_{d-1}=1$ and so it
follows from  Lemmas \ref{lem00} and \ref{remark}, that $G$ is
locally nilpotent.\\ (2) \; if $d=2$ then $c\leq 4$ and so every
group  in the variety $\mathcal{N}_c^{(d)}$ is $4$-Engel. But as
it is mentioned in the proof of Lemma \ref{lem000},  every
$4$-Engel $p$-group is locally nilpotent, where $p\in\{2,3,5\}$.
Now assume that $d\geq 3$, then every group $G$ in the variety
$\mathcal{N}_c^{(d)}$ satisfies the law $V_{d-2}=1$ and so it
follows from Lemmas \ref{lem000} and \ref{remark}, that $G$ is
locally nilpotent. \;\;\;$\Box$

\section{$p$-groups in the variety $\mathcal{N}_c^{(2)}$}

\begin{lem} Let $c\geq 1$ be an integer, $p$ a prime number and  $G$  a finite  $c$-Engel $p$-group.
Suppose that $x,y\in G$ such that $x^{p^n}=y^{p^n}=1$ for some
integer $n>0$. Let $r$ be the integer such that $p^{r-1}<c\leq p^r$. Then:\\
(a) \; if  $p$ is odd and $n>r$ then $[x^{p^{n-1}},y^{p^{n-1}}]=1$.\\
(b) \; if $p=2$ and $n>r+1$ then $[x^{2^{n-1}},y^{2^{n-1}}]=1$. \label{lem1}
\end{lem}
\begin{proof}
Suppose that $K\leq H$ are two normal subgroups of $G$ such that
$H/K$ is elementary abelian and $a$ is an arbitrary element of
$G$. Put $t=aK$ and $V=H/K$. Since $[V, _c t]=1$, we have that
$[V, _{p^r}t]=1$ and $0=(t-1)^{p^r}=t^{p^r}-1$ in $\text{End}(V)$.
Thus $[H, a^{p^r}]\leq K$ for all $a\in G$. Now let
$N=\left<x^{p^r},y^{p^r}\right>$. Then $[H,N]\leq K$ and since
$K,H$ are normal in $G$; $[H,M]\leq K$ where $M=N^G$ the normal
closure of $N$ in $G$. Thus $M$ is a normal subgroup of $G$
centralizing every abelian
 normal section of  $G$. By a result of Shalev \cite{Shalev}; if $p$ is odd then $M$ is powerful and if
$p=2$ then $M^2$ as well as all subgroups of $M^2$ which are normal in $G$, are powerful. Suppose that $p$ is
 odd. Since $M$ is generated by $\{(x^g)^{p^r},(y^g)^{p^r} \;|\; g\in G\}$ and $M$ is powerful,  by
Corollary 1.9 of \cite{LubMan}, $M^{p^{n-r}}$ is generated by $\{(x^g)^{p^n},(y^g)^{p^n} \;|\; g\in G\}$
and so $M^{p^{n-r}}=1$. On the other hand $M^{p^{n-r-1}}$ is powerful by Corollary 1.2 of \cite{LubMan}.
Thus $[M^{p^{n-r-1}},M^{p^{n-r-1}}]\leq (M^{p^{n-r-1}})^p$. Now by Theorem 1.3 of \cite{LubMan},
 we have $(M^{p^{n-r-1}})^p=M^{p^{n-r}}=1$.
Thus $M^{p^{n-r-1}}$ is abelian and the part (a) has been proved. Now assume that $p=2$. As it is mentioned
in the above, since the subgroup $\left<x^{2^{r+1}},y^{2^{r+1}}\right>^G$ of $M$ is normal in $G$, it is also
  powerful and the rest of the proof is similar to the latter case.
\end{proof}
\begin{lem} \label{lem2}
Let $c>1$ be an integer, $p$ a prime number and  $G$  a $p$-group in the variety
 $\mathcal{N}_c^{(2)}$
 and let
 $r$ be the integer
 satisfying $p^{r-1}<c-1\leq p^r$. \\
(a) \; if $p$ is odd then $G^{p^r}$ is locally nilpotent.\\
(b) \; if $p=2$ then $G^{2^{r+1}}$ is locally nilpotent.
\end{lem}
\begin{proof}
Note that $G$ is a $c$-Engel group.
Suppose  $p$ is odd. First we prove that if $x$ is an element of $G$ such that $x^{p^n}=1$ for  some
$n>r$, then  $[a, x^{p^{n-1}},x^{p^{n-1}}]=1$ for all $a\in G$. Let $y=(x^{-1})^a$, then it
 is enough to show that
$[y^{p^{n-1}},x^{p^{n-1}}]=1$. Since $\left<x,a\right>\in\mathcal{N}_c$ then
 $\left<y,x\right>\in\mathcal{N}_{c-1}$; therefore by Lemma \ref{lem1}, we have that
$[y^{p^{n-1}},x^{p^{n-1}}]=1$. Now let $z$ be an arbitrary element of $G$ such that $z^{p^n}=1$
for some positive integer $n$. We prove by induction on $n$ that $z^{p^r}\in HP(G)$ and so it completes
the proof for the case $p$ odd. If $n\leq r$,
 then $z^{p^r}=1\in HP(G)$. Assume that $n>r$ then by induction hypothesis, $z^{p^{r+1}}\in HP(G)$.
 Now by the first part of the proof, $[a,z^{p^r},z^{p^r}]=1$ mod $HP(G)$, for all $a\in G$. So the normal
closure of $z^{p^r}HP(G)$ in  $G/HP(G)$ is abelian. Therefore $\left<z^{p^r}\right>^G$ is  (locally
 nilpotent)-by-abelian, hence Lemma \ref{remark} implies that
 $z^{p^r}\in HP(G)$.\\
The case $p=2$ is similar.
\end{proof}

{\bf Proof of Theorem \ref{thm1}.} Suppose that $p$ is  odd and
every $p$-group of exponent dividing $p^r$ in
$\mathcal{N}_c^{(2)}$ is locally nilpotent. Let $G$ be a $p$-group
in $\mathcal{N}_c^{(2)}$, then  by Lemma \ref{lem2}(a), $G^{p^r}$
is locally nilpotent. By assumption, $G/G^{p^r}$ is locally
nilpotent and since $G$ is a $c$-Engel group, it follows from
Lemma \ref{remark} that $G$ is locally nilpotent. The case $p=2$
is similar and the converse is  obvious.\;\;\;\ $\Box$\\

Now we use this result to some special cases. In fact we prove:
\begin{prop}
Every $2$-group or $3$-group  in the variety $\mathcal{N}_5^{(2)}$ is locally nilpotent.
\end{prop}
\begin{proof}
 By Theorem \ref{thm1}, we must prove that every $2$-group of exponent dividing $8$ and  every $3$-group
of exponent dividing $9$ in  $\mathcal{N}_5^{(2)}$ is  locally nilpotent.
 Suppose that $G$ is a $2$-group of exponent $8$ in  $\mathcal{N}_5^{(2)}$. Let $x,y\in G$,
then $\left<x,y\right>\in \mathcal{N}_5$ and of exponent $8$. It is easy to see that $[x^4,y,x^4,y]=1$. So
by Lemma  \ref{thmhein}, $x^4\in HP(G)$ for all $x\in G$. Therefore $G/HP(G)$ is of exponent
dividing $4$. By a famous result of Sanov (see  14.2.4 of \cite{Robinson}), $G/HP(G)$
is locally nilpotent and so by Lemma \ref{remark}, $G$ is
locally nilpotent.\\
Now suppose that $G$ is  a $3$-group of exponent $9$ in the variety $\mathcal{N}_5^{(2)}$. Let $x,y\in G$,
it is easy to see that $[x^3,y,x^3,y]=1$. The rest of the proof is similar to latter case, but we may use this
well-known result   that every group of exponent $3$ is  nilpotent (see  12.3.5 and 12.3.6 of
\cite{Robinson})
\end{proof}

\end{document}